\documentclass[runningheads,a4paper]{llncs}
\usepackage{graphicx}
\usepackage{amsmath}
\usepackage{bm}

\usepackage{tikz}
\usetikzlibrary{positioning}
\usepackage{textcomp}
\usepackage[hyphens,allowmove]{url}
\usepackage{hyperref}
\usepackage{cleveref}

\begin{document}
\title{Solving System of Nonlinear Equations\\
 with the Genetic Algorithm\\
 and Newton’s Method}
\titlerunning{The Genetic Algorithm and Newton’s Method}
%
\author{Ryuji Koshikawa\inst{1}\thanks{Current affiliation: 
MITA International School, Tokyo, Japan. \url{http://www.mita-is.ed.jp/en/}} \and
Akira Terui\inst{1}\orcidID{0000-0003-0846-3643} \and
Masahiko Mikawa\inst{1}\orcidID{0000-0002-2193-3198}}
\authorrunning{Koshikawa et al.}
%
\institute{%
    University of Tsukuba, Tsukuba, Japan \\
    \email{ryuji0503@math.tsukuba.ac.jp}\\
    \email{terui@math.tsukuba.ac.jp}\\
    \email{mikawa@slis.tsukuba.ac.jp}\\
    \url{https://researchmap.jp/aterui}
}
\maketitle              
\begin{abstract}
    An implementation and an application of the 
    combination of the genetic algorithm and Newton's method for solving a 
    system of nonlinear equations is presented.
    The method first uses the advantage of the robustness of the genetic 
    algorithm for guessing the rough location of the roots, then it uses the 
    advantage of
    a good rate of convergence of Newton’s method. 
    An effective application of the method for the positioning problem of 
    multiple small rovers proposed for the use in asteroid exploration is shown.
\end{abstract}
\section{Introduction}
\label{sec:intro}

Solving a system of nonlinear equations is one of the fundamental problems in science and technology. If the given system consists of algebraic equations, it can be solved with algebraic techniques such as triangularization of the system using Gr\"obner bases \cite{cox-lit-osh2015}. However, if the given system has non-algebraic equations, algebraic methods may not be applicable. 

In such a case, other numerical methods including the genetic algorithm or Newton’s method are used for computing the 
approximate roots (\cite{man-bha-ahm-udd2017}, \cite{rov-val-cas2005}).  
Newton’s method \cite{deu2011} has a good rate of convergence, although the performance of the 
method depends on the initial values, and it may be difficult to find good 
initial values in some cases. On the other hand, a method based on the genetic 
algorithm \cite{mit1998} is sufficiently robust to find global solutions, although the 
convergence rate may not be high, compared with Newton’s method.
Furthermore, a combination of both methods has been proposed 
\cite{kar-wec-fre1998}: 
first by guessing good initial values by the genetic algorithm, then by using Newton’s method for fast convergence. 

In this paper, we present an implementation and an application of the 
combination 
of the genetic algorithm and Newton's method for solving a system of 
nonlinear equations. we demonstrate that the present method finds 
the solutions effectively in the positioning problem of multiple small rovers 
proposed for the use in asteroid exploration. A comparison with the 
root-finding with the genetic algorithm is presented.

The rest of this papers is organized as follows. In \Cref{sec:appl}, 
an application of the present method to the positioning problem of multiple rovers and the formulation of a system of nonlinear equations is presented.
In \Cref{sec:exp}, the result of experiments along with comparison with 
a method using only the genetic algorithm is shown.

\section{Application of the present method to the
positioning problem of multiple rovers to be used in asteroid exploration}
\label{sec:appl}

In this section, we explain an application of the present method to 
the positioning problem of multiple rovers, proposed by one of the present
authors \cite{mik2013}, for asteroid exploration. 

\subsection{Specification of the rovers}
\label{sec:rover-spec}

The rover has
a rectangular shape with approximately 50 mm in length, width, and height,
and multiple rovers are placed on the asteroid for exploration.
Each rover communicate with 
others using radio wave with configuring a wireless mesh network
on a surface of the asteroid. In the radio wave communication,
the rover obtains the received
signal strength indicator (RSSI), which is used for estimating 
relative distances among the rovers.

For communicating with other rovers under good condition even if 
the rover is placed with unusual position (such as upside down), 
the rover has a pair of two antennas on each side, and it uses 
a pair of antennas on the top side for communicating with other
rovers (see \cref{fig:rover-antenna}). 
Let antennas $a$ and $b$ be the pair of the two on the top side, respectively. 

\subsection{Estimating relative distances among the rovers by
using the RSSI}

\begin{figure}[t]
    \begin{minipage}{0.35\textwidth}
        \centering
        \includegraphics[width=0.9\textwidth]{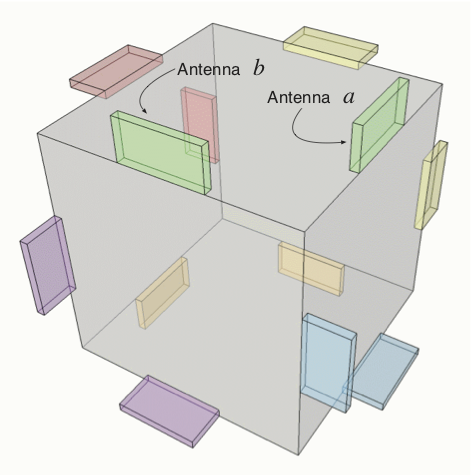}
        \caption{
        The rover equipped with twelve antennas \cite[Fig.\ 3]{mik2013}
        \protect\footnotemark[1].}
        \label{fig:rover-antenna}
    \end{minipage}
    \hfill
    \begin{minipage}{0.55\textwidth}
        \centering
        \includegraphics[width=0.9\textwidth]{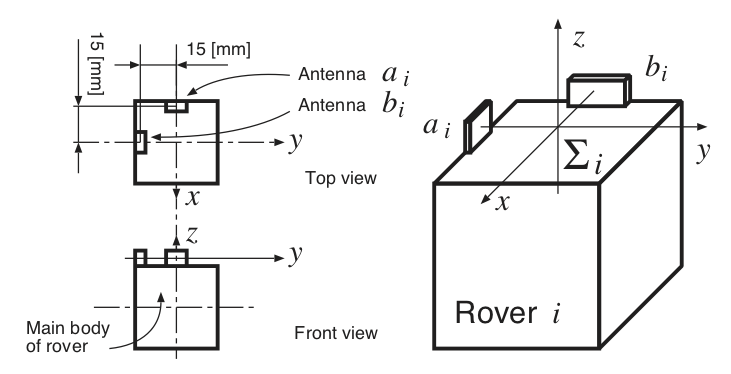}
        \caption{The coordinate system on the rover $i$ \cite[Fig.\ 10]{mik2013}
        \protect\footnotemark[1].}
        \label{fig:antenna-coordinate}
    \end{minipage}
\end{figure}

For estimating relative distances among the rovers by
using the RSSI, for rover $i$ ($i=1,\ldots,n$), 
we define a coordinate system $\sum_i$ as shown in
\cref{fig:antenna-coordinate}.
\footnotetext[1]{\textcopyright\ 2013 IEEE.}

Since each rover communicates with others using antennas on the top side, 
let $a_i$ and $b_i$ be the antennas $a$ and $b$ on the top side, respectively, and let $\sum_{a_i}$ and $\sum_{b_i}$ be coordinate systems with the origin
placed at the center of $a_i$ and $b_i$, respectively. 
Let
\[
    {}^{a_i}\bm{p}_{b_j}=
    \begin{pmatrix}
        {}^{a_i}x_{b_j} \\
        {}^{a_i}y_{b_j} \\
        {}^{a_i}z_{b_j} 
    \end{pmatrix}
    \in R^3
\]
denote the position of the antenna $b_j$ with respect to the coordinate system
$\sum_{a_i}$. Especially, ${}^{1}\bm{p}_{a_i}$ denote the position of the 
antenna $a_i$ with respect to the coordinate system $\sum_{1}$. 
Furthermore, let
${}^1R_{a_i}\in R^{3\times 3}$ denote the orientation of the antenna 
$a_i$ with respect to
the coordinate system $\sum_{1}$.
Then, by using the relationship
\[
    {}^{1}\bm{p}_{a_i} = {}^{1}\bm{p}_{b_j}+{}^{1}R_{b_j} {}^{b_j}\bm{p}_{a_i}, \quad
    {}^{1}\bm{p}_{b_j} = {}^{1}\bm{p}_{a_i}+{}^{1}R_{a_i} {}^{a_i}\bm{p}_{b_j},
\]
we derive the relative position of the antenna $b_j$ with respect to 
the coordinate system $\sum_{a_i}$ and vice versa, as 
\[
    {}^{b_j}\bm{p}_{a_i}=({}^{1}\bm{p}_{a_i}-{}^{1}\bm{p}_{b_{j}})
    {}^{b_j}R_{1},\quad 
    {}^{a_i}\bm{p}_{b_j}=({}^{1}\bm{p}_{b_{j}}-{}^{1}\bm{p}_{a_{i}})
    {}^{a_i}R_{1},
\]
noting that ${}^{b_j}R_{1}=({}^{1}R_{b_j})^{-1}$ and 
${}^{a_i}R_{1}=({}^{1}R_{a_i})^{-1}$.

\begin{figure}[t]
    \centering
    \includegraphics{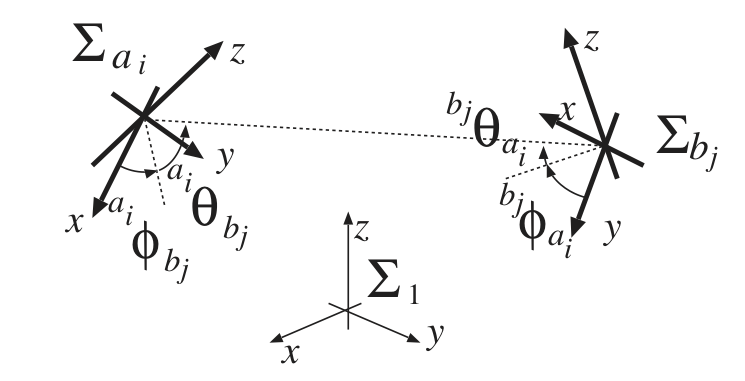}
    \caption{The horizontal and the elevation angles of rovers
    \cite[Fig.\ 11]{mik2013}
    \protect\footnotemark[2].}
    \label{fig:coordinate-systems}
\end{figure}

Next, let $^{a_i}\phi_{b_j}$ and $^{a_i}\theta_{b_j}$ denote
the horizontal and the elevation angle of the antenna $b_j$
with respect to the coordinate system $\sum_{a_i}$, respectively, 
as shown in \cref{fig:coordinate-systems}. 
\footnotetext[2]{\textcopyright\ 2013 IEEE.}
Then, by using the relative position 
${}^{a_i}\bm{p}_{b_j}={}^t({}^{a_i}x_{b_j},{}^{a_i}y_{b_j},
{}^{a_i}z_{b_j} )$,
$^{a_i}\phi_{b_j}$ and $^{a_i}\theta_{b_j}$ are derived as 
\[
    {}^{a_i}\phi_{b_j}=\arctan
    \left(
        \frac{{}^{a_i}y_{b_j}}{{}^{a_i}x_{b_j}}
    \right),
    \quad
    {}^{a_i}\theta_{b_j}=\arctan
    \left(
        \frac{{}^{a_i}z_{b_j}}{\sqrt{({}^{a_i}x_{b_j})^2
        +({}^{a_i}y_{b_j})^2}}
    \right).
\]
The horizontal and the elevation angle of the antenna $a_i$
with respect to the coordinate system $\sum_{b_j}$, respectively,
are derived similarly.

Let $r_{a_i\_b_j}$ be the RSSI between the antenna $a_i$ and $b_j$. 
By the discussion above, the mathematical model of $r_{a_i\_b_j}$ are
expressed as 
\begin{multline}
    \label{eq:raibj}
    r_{a_i\_b_j}({}^{a_i}x_{b_j},{}^{a_i}y_{b_j},{}^{a_i}z_{b_j},
    {}^{a_i}\phi_{b_j},{}^{a_i}\theta_{b_j},
    {}^{b_j}\phi_{a_i},{}^{b_j}\theta_{a_i}) \\
    =r^{\prime}_{a_i,b_j}{}(^{a_i}x_{b_j},{}^{a_i}y_{b_j},{}^{a_i}z_{b_j}) 
    +r_h({}^{a_i}\phi_{b_j})+r_h({}^{b_j}\phi_{a_i})
    +r_v({}^{a_i}\theta_{b_j})+r_v({}^{b_j}\theta_{a_i}),
\end{multline}
where 
\begin{align*}
    r'(x,y,z)&=-14.69\log_{10}(\sqrt{x^{2}+y^{2}+z^{2}}+0.31)-49.17, \\
    r_h(\phi)&=\frac{5}{2}(\cos (2\phi)-1), \\
    r_v(\theta)&=25
        \left(
            \frac{\cos(\frac{5}{2}\pi)\cos(\frac{5}{2}\pi-|\theta|)}{\sin (\frac{5}{2}\pi-|\theta|)}
        \right)-1.
\end{align*}
We find the values of 
${}^{a_i}x_{b_j},{}^{a_i}y_{b_j},{}^{a_i}z_{b_j},
{}^{a_i}\phi_{b_j},{}^{a_i}\theta_{b_j},{}^{b_j}\phi_{a_i},{}^{b_j}\theta_{a_i}$
by solving \cref{eq:raibj} with respect to these variables. 
In the original article \cite{mik2013}, the author solves the equation by using 
the genetic algorithm in whole using GAlib \cite{galib}. 
In the genetic algorithm, the evaluation function is defined as follows. 
Let $\bar{r}_{a_i\_b_j}$ and $\hat{r}_{a_i\_b_j}$ be the measured and the 
estimated RSSI values, respectively. Then, the evaluation function $f(r)$ 
is defined as 
\begin{multline}
    f(r)=\sum^{n}_{i=1} \sum^{n}_{j=i+1} 
    \left\{
        (\bar{r}_{a_i\_a_j}-\hat{r}_{a_i\_a_j})^2
        +(\bar{r}_{a_i\_b_j}-\hat{r}_{a_i\_b_j})^2
    \right. 
        \\
    \left.
        +(\bar{r}_{b_i\_a_j}-\hat{r}_{b_i\_a_j})^2
        +(\bar{r}_{b_i\_b_j}-\hat{r}_{b_i\_b_j})^2 
    \right\}.
    \label{eq:evalf}
\end{multline}

\subsection{Setting up nonlinear equations 
with Newton's method and the genetic algorithm}
\label{sec:settingup-equations}

In this paper, let us assume that 
all the rovers have the same direction and they are placed on the 
$xy$ plane so that the elevation angle ${}^{a_i}\theta_{b_j}$ are 
omitted. 

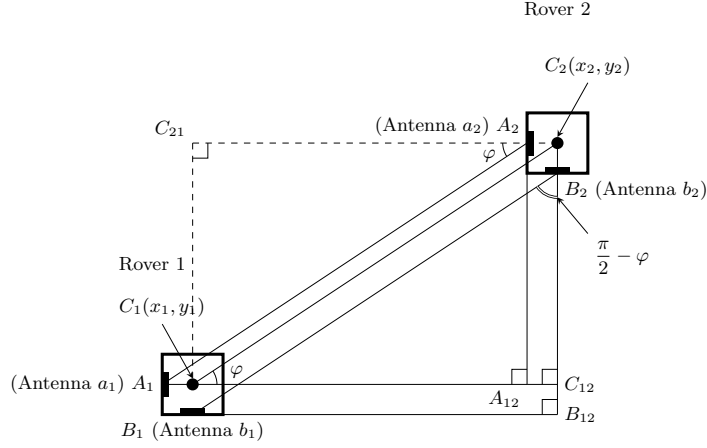
\begin{figure}[t]
    \centering
    \scalebox{0.8}{
    \begin{tikzpicture}[>=stealth]
  
      \coordinate[label = left:{Rover 1}] () at (0.5, 2.5);
      \draw[line width = 0.5mm] (0,0) rectangle (1,1);
      \coordinate[label = left:{(Antenna $a_1$) $A_1$}] (a1) at (0,0.5);
      \coordinate[label = below:{$B_1$ (Antenna $b_1$)}] (b1) at (0.5,0);
      \coordinate (c1) at (0.5,0.5);
      \draw[->] (0,1.5) -- (0.5,0.6);
      \coordinate[label = above:{$C_1(x_1,y_1)$}] (c1l) at (0,1.5);
      
      \fill[black] (0,0.3) rectangle (0.1,0.7);
      \fill[black] (0.5, 0.5) circle [radius = 1mm];
      \fill[black] (0.3,0) rectangle (0.7,0.1);
      
      \coordinate[label = above:{Rover 2}] () at (6.5, 6.5);
      \draw[line width = 0.5mm] (6,4) rectangle (7,5);
      \coordinate[label = 135:{(Antenna $a_2$) $A_2$}] (a2) at (6,4.5);
      \coordinate[label = 315:{$B_2$ (Antenna $b_2$)}] (b2) at (6.5,4);
      \coordinate (c2) at (6.5,4.5);
  
      \coordinate[label = 250:{$A_{12}$}] (a12) at (6,0.5);
      \coordinate[label = right:{$B_{12}$}] (b12) at (6.5,0);
      \coordinate[label = right:{$C_{12}$}] (c12) at (6.5,0.5);
      \draw[->] (7,5.5) -- (6.5,4.6);
      \coordinate[label = above:{$C_2(x_2,y_2)$}] (c2l) at (7,5.5);
      
      \fill[black] (6.5, 4.5) circle [radius = 1mm];
      \fill[black] (6,4.3) rectangle (6.1,4.7);
      \fill[black] (6.3,4) rectangle (6.7,4.1);
     
      \coordinate[label = 135:{$C_{21}$}] (c21) at (0.5,4.5);

      \draw (a1) -- (a2);
      \draw (a1) -- (a12);
      \draw (a12) -- (c12);
      \draw (a12) -- (a2);
      
      \draw (b1) -- (b2);
      \draw (b1) -- (b12);
      \draw (b12) -- (b2);
      \draw [double] (6.5,3.6) arc [start angle = 270, end angle = 213.7, radius = 4mm];
      \draw[->] (7,3) -- (6.5,3.6);
      \coordinate[label = -30:{$\displaystyle{\frac{\pi}{2}-\varphi}$}] () at (7,3);
      
      \draw [dashed] (c21) -- (c2);
      \draw [dashed] (c21) -- (c1);
      \draw (0.75,4.5) -- (0.75,4.25) -- (0.5,4.25);

      \draw (5.6,4.5) arc [start angle = 180, end angle = 213.7, radius = 4mm];
      \coordinate[label = 215:{$\varphi$}] () at (5.6,4.5);
     
      \draw (c1) -- (c2);
      \draw (b2) -- (c2);
      \draw (0.9,0.5) arc [start angle = 0, end angle = 33.7, radius = 4mm];
      \coordinate[label = 45:{$\varphi$}] () at (1,0.5);
      
      \draw (6.25,0) -- (6.25,0.25) -- (6.5,0.25);
      \draw (6.25,0.5) -- (6.25,0.75) -- (6.5,0.75);
      \draw (5.75,0.5) -- (5.75,0.75) -- (6,0.75);
  
    \end{tikzpicture}
    }
    \caption{Positions of two rovers.}
    \label{fig:rover12}
  \end{figure}
Nonlinear equations for estimating relative distance 
and the rotation angle between two rovers are derived, as shown in 
\cref{fig:rover12}. Let $\varphi$ be the acute angle between lines 
$C_1C_{12}$ and $C_1C_2$. Note that triangles $A_1A_2A_{12}$, 
$B_1B_2B_{12}$ and $C_1C_2C_{12}$ are congruent. 
Since $\angle C_2C_1C_{12}=\angle A_2A_1A_{12}$ and 
$\angle A_2A_1A_{12}$ and $\angle A_1A_2C_{21}$ are the complex angle,
we have $\angle A_1A_2C_{21}=\varphi$. Furthermore, 
$\angle B_2B_1B_{12}=\varphi$ shows that 
$\angle B_1B_2B_{12}=\frac{\pi}{2}-\varphi$.
Thus, the RSSI is expressed as 
\begin{align}
    \label{eq:ra1a2}
    & r_{a_1\_{a_2}}({}^{a_1}x_{a_2},{}^{a_1}y_{a_2},{}^{a_1}\varphi_{a_2},{}^{a_2}\varphi_{a_1}) \nonumber \\
    &=-14.69\log_{10}(\sqrt{({}^{a_1}x_{a_2})^2+({}^{a_1}y_{a_2})^2}+0.31)-49.17+5(\cos 2\varphi -1), \\
    \label{eq:rb1b2}
    & r_{b_1\_{b_2}}({}^{b_1}x_{b_2},{}^{b_1}y_{b_2},{}^{b_1}\varphi_{b_2},{}^{b_2}\phi_{b_1}) \nonumber \\
    &=-14.69\log_{10}(\sqrt{({}^{b_1}x_{b_2})^2+({}^{b_1}y_{b_2})^2}+0.31)-49.17+5(\cos 2(\frac{\pi}{2}-\varphi) -1).
\end{align}
By using $A_1A_2=B_1B_2$ and rewriting $\cos 2(\frac{\pi}{2}-\varphi)$ as
\begin{align*}
    \cos\left(2\left(\displaystyle{\frac{\pi}{2}-\varphi}\right)\right)
    &= \cos^2\left(\displaystyle{\frac{\pi}{2}-\varphi}\right)
    - \sin^2\left(\displaystyle{\frac{\pi}{2}-\varphi}\right) \\
    &= \sin^2\varphi-\cos^2\varphi
    = -\cos(2\varphi),
\end{align*}
we have
\begin{equation}
    \label{eq:varphi}
    r_{a_1\_a_2}-r_{b_1\_b_2}=10\cos(2\varphi),
\end{equation}
which gives the value of $\varphi$. Then, by solving a system of 
nonlinear equations
\begin{equation}
    \label{eq:position}
    \begin{split}
        r_{a_1\_a_2}({}^{a_1}x_{a_2},{}^{a_1}y_{a_2},\varphi) &=-14.69 
        \log_{10} (\sqrt{({}^{a_{1}}x_{a_{2}})^{2}+
        ({}^{a_{1}}y_{a_{2}})^{2}}+0.31)
        \\
        &\quad -49.17+5 (\cos 2\varphi -1), \\
        \varphi &=\arctan\left(\dfrac{{}^{a_1}y_{a_2}}{{}^{a_1}x_{a_2}}\right),
    \end{split}
\end{equation}
we find the relative position of $(x_2,y_2)$ 
with respect to $(x_1,y_1)$.

In this paper, we assume that there exists a rover 
on the origin 
and we estimate the position of the other rovers by computing the distance 
and the orientation angle of those from the one on the origin.
In the genetic algorithm, the evaluation function 
is defined as in \cref{eq:evalf}. In Newton’s method
applied to \cref{eq:position}, the Jacobian matrix becomes as 
\begin{multline}
    J=
    \left(
    \begin{array}{c}
        \dfrac{-14.69 ({}^{a_1}x_{a_2}) \ln 10}{({}^{a_1}x_{a_2})^2
        +({}^{a_1}y_{a_2})^2+0.31\sqrt{({}^{a_1}x_{a_2})^2+({}^{a_1}y_{a_2})^2}} 
        \\
        \dfrac{{}^{a_1}x_{a_2}}{({}^{a_1}x_{a_2})^2+({}^{a_1}y_{a_2})^2}
    \end{array}
    \right.
    \\
    \left.
    \begin{array}{c}
        \dfrac{-14.69 ({}^{a_1}y_{a_2}) \ln 10}{({}^{a_1}x_{a_2})^2
        +({}^{a_1}y_{a_2})^2
        +0.31\sqrt{({}^{a_1}x_{a_2})^2+({}^{a_1}y_{a_2})^2}}  \\
        \dfrac{{}^{a_1}x_{a_2}}{({}^{a_1}x_{a_2})^2+({}^{a_1}y_{a_2})^2} 
    \end{array}
    \right).
\end{multline}

\section{Experiments}
\label{sec:exp}

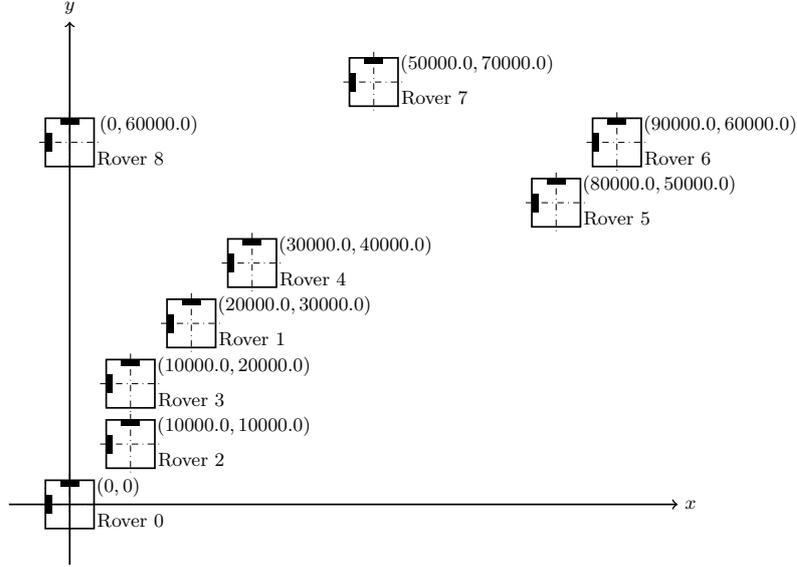
\begin{figure}[t]
    \centering
    \scalebox{0.8}{
        \begin{tikzpicture}[%
                rover/.pic = {
                    \node[rectangle, draw, thick, minimum size=0.8cm] (rover0) {} (0,0);
                    \fill[black] (-0.4,-0.15) rectangle (-0.3,0.15);
                    \fill[black] (-0.15,0.3) rectangle (0.15,0.4);
                    \draw[dash dot] (0,0.5)--(0,-0.5); 
                    \draw[dash dot] (-0.5,0)--(0.5,0); 
                }
            ]

            \draw[->,thick] (-1,0)--(10,0) node[right]{$x$};
            \draw[->,thick] (0,-1)--(0,8) node[above]{$y$};
            
            \draw (0,0) pic {rover} 
                node [label={[xshift=1cm, yshift=-0.6cm]Rover 0}] {}
                node [label={[xshift=0.8cm, yshift=-0.1cm] $(0,0)$}] {}
                ;
            \draw (2,3) pic {rover} 
                node [label={[xshift=1cm, yshift=-0.6cm]Rover 1}] {}
                node [label={[xshift=1.7cm, yshift=-0.1cm] $(20000.0,30000.0)$}] {};
            \draw (1,1) pic {rover} 
                node [label={[xshift=1cm, yshift=-0.6cm]Rover 2}] {}
                node [label={[xshift=1.7cm, yshift=-0.1cm] $(10000.0,10000.0)$}] {};
            \draw (1,2) pic {rover} 
                node [label={[xshift=1cm, yshift=-0.6cm]Rover 3}] {}
                node [label={[xshift=1.7cm, yshift=-0.1cm] $(10000.0,20000.0)$}] {};              
            \draw (3,4) pic {rover} 
                node [label={[xshift=1cm, yshift=-0.6cm]Rover 4}] {}
                node [label={[xshift=1.7cm, yshift=-0.1cm] $(30000.0,40000.0)$}] {}; 
            \draw (8,5) pic {rover} 
                node [label={[xshift=1cm, yshift=-0.6cm]Rover 5}] {}
                node [label={[xshift=1.7cm, yshift=-0.1cm] $(80000.0,50000.0)$}] {};
            \draw (9,6) pic {rover} 
                node [label={[xshift=1cm, yshift=-0.6cm]Rover 6}] {}
                node [label={[xshift=1.7cm, yshift=-0.1cm] $(90000.0,60000.0)$}] {};  
            \draw (5,7) pic {rover} 
                node [label={[xshift=1cm, yshift=-0.6cm]Rover 7}] {}
                node [label={[xshift=1.7cm, yshift=-0.1cm] $(50000.0,70000.0)$}] {}
                ; 
            \draw (0,6) pic {rover} 
                node [label={[xshift=1cm, yshift=-0.6cm]Rover 8}] {}
                node [label={[xshift=1.3cm, yshift=-0.1cm] $(0,60000.0)$}] {}
                ; 
        \end{tikzpicture}
    }
    \caption{Positions of the rovers (unit: mm).}
    \label{fig:rover-experiment}
\end{figure}

We have tested the method for estimating positions of 
rovers placed on the $xy$ plane as shown in \cref{fig:rover-experiment}.
For the implementation of the method, GAlib \cite{galib} was 
used for the genetic algorithm and our implementation was used
for Newton’s method.

Position of Rover $i$ ($i=1,\ldots,8$) is estimated by computing
the distance and the horizontal angle between Rover $0$ placed
on the origin. 
In the experiments, the genetic algorithm was executed with the
following settings: 
the number of population was set to 100, 
the number of generation was set to 200,
single-point crossover and roulette wheel selection were used
with the crossover rate set to 0.9,
and the mutation rate was set to 0.01.  
After evolving the population for repeating the computation for 
$10\times 10=100$ times (which is denoted by $10\times 10$), 
the best individuals which minimize $f(r)$ were
used as the initial values
for Newton's method. Newton's method was terminated either when
the magnitude of the updated value became less than $1.0\times 10^{-10}$
or the number of iterations reached 100.

The computing environment is as follows: 
Intel Xeon E5607 2.27GHz, RAM 48 GB, Linux 3.16.0, GCC 9.2.1.

\subsection{Experiment 1: estimation of rovers with the present method}
\label{sec:exp-ga-newton}
\Cref{tab:ga-10x10} shows the result of the genetic algorithm computation. 
The columns are as follows: `$i$' is the index of the rover, `Actual position'
and `Estimated position' is the actual and the estimated position 
of Rover $i$, respectively, `Relative error' is 
\begin{equation}
    \label{eq:relative-error}
    |d_i-d'_i|/d_i,
\end{equation}
where $d_i$ and $d'_i$ is the actual and the estimated distance of Rover $i$ 
from the origin, respectively. 
The bottom row (with `Average' in the column of `Actual position') shows 
the average of the realtive errors.
Note that the estimated position is used as the initial 
value of Newton's method.
\begin{table}[t]
    \caption{The result of estimation by the genetic algorithm
     with the setting of $10\times 10$.}
    \label{tab:ga-10x10}
    \begin{center}
        \begin{tabular}{l|l|l|l}
            \hline
            $i$ & Actual position & 
            \begin{tabular}{l}
                Estimated position\\
                (The initial value for Newton's method)
            \end{tabular}
            & Relative error \eqref{eq:relative-error}
            \\
            \hline
            1 & $(20000.0,30000.0)$ & $(18411.987305,33178.985596)$ &
            $0.052413$
            \\
            2 & $(10000.0,10000.0)$ & $(10122.98584,10196.990967)$ &
            $0.016006$
            \\
            3 & $(10000.0,20000.0)$ & $(11677.993774,20789.993286)$ &
            $0.066396$
            \\
            4 & $(30000.0,40000.0)$ & $(26173.995972,43239.990234)$ &
            $0.010896$ 
            \\
            5 & $(80000.0,50000.0)$ & $(77340.995789,55555.999756)$ &
            $0.009400$
            \\
            6 & $(90000.0,60000.0)$ & $(89379.997253,64320.999146)$ &
            $0.018041$
            \\
            7 & $(50000.0,70000.0)$ & $(50185.989380,75285.995483)$ &
            $0.051808$
            \\
            8 & $(0,60000.0)$ & $(2766.006470,64238.998413)$ &
            $0.071642$ 
            \\
            \hline
            --- & Average & --- & $0.037075$
            \\
            \hline
        \end{tabular}
    \end{center}
\end{table}

\Cref{tab:newton} shows the result of Newton's method after 
the genetic algorithm. The columns `$i$', `Actual position', 
`Estimated position' and `Relative error' are the same as those
in \Cref{tab:ga-10x10}. The column `\#iterations' is 
the number of iterations executed for estimating the position. 
For Rover 5, Newton's method has terminated after the number of iteration 
becomes 100, while, for the other rovers, it has terminated
the magnitude of the updated value has become less than $1.0\times 10^{-10}$.
We see that for all the rovers except for Rover 5, the relative error 
decreases after applying Newton's method. Especially, for Rovers 
2, 3 and 8, the relative error decreases significantly.
\begin{table}[t]
    \caption{The result of estimation by Newton's method.}
    \label{tab:newton}
    \begin{center}
        \begin{tabular}{l|l|l|l|r}
            \hline
            $i$ & Actual position & Estimated position &  Relative error
            \eqref{eq:relative-error} & 
            \#iterations
            \\
            \hline
            1 & $(20000.0,30000.0)$ & $(19511.732516,29804.663731)$ & 0.011985
            & 99 \\
            2 & $(10000.0,10000.0)$ & $(10000.000000,10000.000000)$ & 
            $<1.0\times10^{-11}$ & 95 \\
            3 & $(10000.0,20000.0)$ & $(10000.000000,20000.000000)$ &
            $<1.0\times10^{-11}$ & 98 \\
            4 & $(30000.0,40000.0)$ & $(29120.352787,39684.351546)$ &
            0.015552 & 98 \\
            5 & $(80000.0,50000.0)$ & $(76602.187691,50147.903356)$ &
            0.029496 & 100 \\
            6 & $(90000.0,60000.0)$ & $(91596.396916,59963.917465)$ &
            0.01213 & 88 \\
            7 & $(50000.0,70000.0)$ & $(51871.537800,70688.990493)$ &
            0.019247 & 96 \\
            8 & $(0,60000.0)$ & $(0.0,60000.000000)$ &
            $<1.0\times10^{-11}$ & 88 \\
            \hline
            --- & Average & --- & $0.011051$ & ---
            \\
            \hline
        \end{tabular}
    \end{center}
\end{table}

\subsection{Experiment 2: estimation of rovers using only the genetic algorithm
and discussions}

For comparison, we have conducted another experiment for estimating 
the position of the rovers using only the genetic algorithm.
In this experiment, after executing the genetic algorithm for 
$n\times n=n^2$ times (which is denoted by $n\times n$),
the best population which minimizes $f(r)$ in \cref{eq:evalf} was
chosen.
We have solved the equations by this method for $n=20,30,40,50,60,70,80$.
All other settings of the genetic algorithm are the same as those 
in the previous experiment. 
\Cref{tab:ga-80x80} shows the average of relative errors 
and the computing time of the experiment.
Note that, for comparison,  the result for $n=10$ is taken 
from \Cref{tab:ga-10x10} and the bottom row shows the result of Experiment 1.

Comparison of the average of relative errors of the estimated positions of the 
rovers in both experiments
(\Cref{tab:ga-80x80}) shows that
the present method can estimate the position of the rovers with 
better accuracy on average than  
the method using only the genetic algorithm.
Furthermore, we see that the present method
is significantly more efficient than the method using only the genetic
algorithm. Thus, we conclude that the method of combining the genetic
algorithm and Newton's method estimates the position 
of rovers more efficiently than methods using only the genetic algorithm,
with better accuracy on average.
\begin{table}[t]
    \caption{The result of estimation by the genetic algorithm
     with the setting of $n\times n$, comparison with the result of 
     the method of combining the genetic algorithm and Newton's method.}
    \label{tab:ga-80x80}
    \begin{center}
        \begin{tabular}{c|c|c}
            \hline
            $n$ & The average of relative errors & Computing time (sec.) 
            \\
            \hline
            10 & $0.037075$ & $2516.31199$ 
            \\
            20 & $0.024174$ & $10047.422038$ 
            \\
            30 & $0.024484$ & $8614.607373$ 
            \\
            40 & $0.020342$ & $15288.596504$ 
            \\
            50 & $0.024484$ & $23874.125484$ 
            \\
            60 & $0.024605$ & $34370.2462$ 
            \\
            70 & $0.032624$ & $46743.015786$ 
            \\
            80 & $0.016325$ & $61102.5245$ 
            \\
            \hline
            GA$+$Newton & $0.011051$ & $2516.31389$ 
            \\
            \hline
        \end{tabular}
    \end{center}
\end{table}

\section{Concluding remarks}
\label{sec:remark}

In this paper, we have demonstrated a numerical method combining the genetic 
algorithm and Newton's method for solving a system of nonlinear equations.
Our method uses the genetic algorithm for a global search of approximate 
roots, then it uses Newton's method for fast convergence. 
We have formulated the problem of estimating the position of rovers used in
asteroid exploration into a system of nonlinear equations for the use of 
the present method.
The experiments 
have shown that the present method computes the roots with almost the same
accuracy and significantly better efficiency than a method using only the 
genetic algorithm.

Future research direction includes the following.
\begin{enumerate}
    \item Since our formulation of a system of equations \eqref{eq:ra1a2} uses 
    arctangent function, we can estimate the position of rovers located only in 
    the first and the 4th quadrant. Another formulation of equations for estimating the position of rovers located in the second and the third quadrant will be needed.
    \item In this paper, we assumed that the rover is placed horizontally
    on the $xy$ plane. However, in the actual exploration, the surface
    of the asteroid is probably not flat, a formulation of the estimation 
    problem in the $xyz$ space with the elevation angle will be needed.
    \item In our method, we first estimated the orientation angle of 
    the rovers from the RSSI values as in \eqref{eq:varphi}, then 
    solved \cref{eq:ra1a2}. Since the 
    error in the measured RSSI values affects the coefficients in 
    \cref{eq:ra1a2}, estimation of errors in the roots of \cref{eq:ra1a2}
    caused by the error in the measured RSSI values will be needed.
\end{enumerate}

\bibliographystyle{splncs04}
\bibliography{icms-2020-newton-genetic-paper}
\end{document}